\documentclass{amsart}
\usepackage{amssymb,latexsym}
\usepackage{enumerate}

\newtheorem{thm}{Theorem}[section]
\newtheorem{cor}[thm]{Corollary}
\newtheorem{prp}[thm]{Proposition}

\theoremstyle{definition}

\newtheorem{rem}[thm]{Remark}

\begin{document}

\title[generalized compositions]{Generalized compositions with a fixed number of parts}

\author[M. Janji\'c]{Milan Janji\'c}
\address{Departments for Mathematics and Informatics\\ University of Banja Luka\\
51000 Banja Luka, Republic of Srpska, BA.}
\email{agnus@blic.net}

\date{}

\begin{abstract}
  We investigate compositions of a positive integer with a fixed number of parts, when there are several  types of each natural number. These compositions produce new  relationships among binomial coefficients, Catalan numbers, and numbers of the Catalan triangle.

\end{abstract}



\maketitle

\section{Introduction}
 A $k$-tuple  $(i_1,i_2,\ldots,i_k)$ of positive integers, such that $i_1+i_2+\cdots+i_k=n,$ is called a
composition of $n$ with $k$ parts.  In \cite{mil}, the following generalization of compositions is considered:
Let $\mathbf b=(b_1,b_2,\ldots,)$ be a  sequence of nonnegative integers, and let $n$ be a positive integer.
The composition of  $n$ is a $k$-tuple $(i_1,i_2,\ldots,i_k)$ such that $i_1+i_2+\cdots+i_k=n,$ assuming that there are $b_1$ different types of  $1$, $b_2$ different types of  $2$, and so on. We call such a composition the generalized composition of $n$ with $k$ parts.

The generalized compositions extend several types of compositions which are investigated in some earlier papers.
First of all this is the case with  usual compositions, which are obtained when $b_i=1$ for each $i.$
In  \cite{deu}, the author considers  the compositions in which there are two different types of $1,$ and one type of each other  natural number. Next, in \cite{ag}, the case $b_i=i,\;(i=1,2,\ldots)$ is investigated.

The generalized compositions may  be  described as the colored compositions, in which the  part $i$ is colored by one of  $b_i$  colors. Different kinds of compositions have already been called  colored compositions. For example,   the $m$-colored compositions, as they are  defined in \cite{drk}, are, freely speaking,  the generalized compositions in which $b_i\in\{\omega,\omega^2,\ldots \omega^{m-1}\},$ where $\omega$ is a primitive $m$th root of $1.$ As well, the composition in which $b_i=i,$ for any $i$, considered in \cite{ag}, is also called an $m$-colored compositions.
The above-mentioned compositions, as well as many other interesting results  on compositions can be found in a recently-published book \cite{hu}.

 In \cite{mil},  several recursions and some closed  formulas for the number of all generalized compositions are obtained.

  In this paper, we investigate the generalized compositions with a fixed number of parts.
The paper is organized as follows. In Section 2 we  outline some basic properties of the generalized compositions with a fixed number of parts. We also show that they extend the notion of the matrix composition, considered in \cite{mu}.
Then we derive several recurrence equations  and closed formulas, by choosing for $b_i$ different functions of $i.$ In particular, we obtain the formula for  the number of $n$-colored compositions, given in \cite{drk}, as well as the formula for the number of $n$-colored compositions, given in \cite{ag}.
Section 3  deals with  the case when $b_i$ is a binomial coefficient. Several closed formulas will be derived. Also, if $b_i$ is of the form $b_i={i+p-1\choose q},$ we prove that the numbers of all generalized compositions satisfy a homogenous recurrence equation with constant coefficients, of order $q+1.$
In particular, the $m$-matrix compositions satisfy such a recurrence equation. For the case $p=1,$ we derive a closed formula for
both the number of the generalized compositions with a fixed number of parts and for the number of all generalized compositions.
In Section 4, we investigate relationships of the generalized compositions with the Catalan numbers. Finally, a result which connects the Catalan numbers, the numbers of the Catalan triangle, and the binomial coefficients is derived.

\section{Some preliminary results}
Let $\mathbf b=(b_0,b_1,\ldots)$ be a sequence of nonnegative integers, and $n,k$ be positive integers. We let  $C^{(\mathbf b)}(n,k)$ denote the number of the generalized compositions of $n$ with $k$ parts.
We also define $C^{(\mathbf b)}(0,0)=1,\;C^{(\mathbf b)}(i,0)=0,\;(i>0).$

In \cite{mil}, the number of all generalized compositions of $n$ is denoted by $C^{(\mathbf b)}(n).$ Obviously,
\begin{equation}\label{sve}C^{(\mathbf b)}(n)=\sum_{k=1}^nC^{(\mathbf b)}(n,k).\end{equation}

In the following two propositions we state some basic properties of  the generalized compositions.
\begin{prp}\label{pp1} The following equations are true:
   \[C^{(\mathbf b)}(i,1)=b_i,\;(i=1,2,\ldots),\;C^{(\mathbf b)}(n,n)=b_1^n,\;C^{(\mathbf b)}(n,k)=0,\;(k>n).\]
\end{prp}
\begin{proof} All equations are easy to verify.
\end{proof}
\begin{prp}\label{rec} The following recursions are true:
\begin{equation}\label{r1}C^{(\mathbf b)}(n,k)=\sum_{i=1}^{n-k+1}b_{i}C^{(\mathbf b)}(n-i,k-1).\;(k\leq n).\end{equation}
\begin{equation}\label{r2}C^{(\mathbf b)}(n)=\sum_{i=1}^{n}b_{i}C^{(\mathbf b)}(n-i),\end{equation}
providing that $C^{(\mathbf b)}(0)=1.$
\end{prp}
\begin{proof} Equation (\ref{r1}) is true since there are $b_iC^{(\mathbf b)}(n-i,k-1)$ generalized compositions ending with one of the $i$'s, for $i=1,\ldots,n-k+1.$ A similar argument proves equation (\ref{r2}).
\end{proof}
We next prove that the matrix compositions, considered in \cite{mu}, are a particular case of the generalized compositions.
 A $k$- matrix composition of $n$ is a matrix with $k$ rows, which entries are nonnegative integers,   no column consists of zeroes only, and the sum of all entries  equals $n.$
We let  $MC(n)$ denote its number.
\begin{prp}\label{pxx} If $b_i={i+k-1\choose i},\;(1=1,2,\ldots),$ then
\[MC(n)=C^{(\mathbf b)}(n).\]
\end{prp}
\begin{proof}
It is a well-known that, for a given positive integer  $k,$ the equation $x_1+x_2+\cdots+x_k=i$ has ${i+k-1\choose i}$ nonnegative solutions. This means that ${i+k-1\choose k}MC(n-i)$ is the number of  $k$-compositions of $n,$ ending with a column in which the sum of all elements equals $i.$
Taking $MC(0)=1$ we obtain
\[MC(n)=\sum_{i=1}^n{i+k-1\choose i}MC(n-i),\]
Comparing this equation with (\ref{r2}) we easily conclude that
\[MC(n)=C^{(\mathbf b)}(n),\] and the proposition is proved.
\end{proof}

In the rest of this section we shall choose for $b_i$ different functions of $i$ and obtain several closed formulas.
We first consider the case when $\mathbf b$ is a constant sequence.
\begin{prp}\label{p21} Let $n,p$ be positive integers, and let $b_i=p,\;(i=1,2,\ldots).$ Then
\[C^{(\mathbf b)}(n,k)=p^k{n-1\choose k-1}.\]
\end{prp}
\begin{proof} In this case, the connection between
compositions and generalized compositions is simple. From a composition of $n$ with $k$ parts we obtain $p^k$ different generalized compositions with $k$ parts, since each part my take $p$ different values. In this way we obtain all generalized compositions of $n$ with $k$ parts.  Moreover,
there are ${n-1\choose k-1}$ compositions of $n$ with $k$
parts, and the proposition follows.
\end{proof}

\begin{cor} In the conditions of Proposition \ref{p21} we have
\[C^{(\mathbf b)}(n)=p(1+p)^{n-1}.\]
\end{cor}
\begin{proof}
The formula (\ref{sve}) now takes the form:
\[C^{(\mathbf b)}(n)=\sum_{k=1}^n{n-1\choose k-1}p^k,\]
and the assertion follows from the binomial formula.
\end{proof}
\begin{rem} The number $C^{(\mathbf b)}(n)$ from the preceding corollary  equals the number of the $p$-colored compositions, as they are defined in \cite{drk}.
\end{rem}
Next, we investigate the case when $\mathbf b$ is a constant sequence with several leading zeroes.
\begin{prp}\label{p22} Let $p,m,n$ be   positive integers, and let  $b_i=0,\;(i=1,2,\ldots,m-1),\;b_i=p,\;(i\geq m).$
Then
\[C^{(\mathbf b)}(n,k)=p^k{n-(m-1)k-1\choose k-1}.\]
\end{prp}
\begin{proof} In this case, we consider  the set $X$ of the generalized compositions of $n$ with $k$ parts, all of which are $\geq m.$
There is a bijection between the set $X$ and  the set $Y$ of  the generalized compositions of $n-(m-1)k$ with $k$ parts, which are considered in Proposition \ref{p21}. Namely, subtracting $m-1$ from each term of an  element of $X,$  we obtain an element of $Y.$  Conversely,  adding $m-1$ to each term of an arbitrary element of $Y,$ we obtain an element of $X.$
The proposition now follows from Proposition \ref{p21}.
\end{proof}
As an immediate consequence of (\ref{sve}) we state
\begin{cor} In the conditions of Proposition \ref{p22} we have
\[C^{(\mathbf b)}(n)=\sum_{k=1}^n{n-(m-1)k-1\choose k-1}p^k.\]
\end{cor}

We shall now consider the case when $b_i$ is an exponential function of $i.$
\begin{prp}\label{p23} Let $p,n,k$ be  positive integers,  and let $b_i=p^{i-1},\;(i=1,2,\ldots).$  Then
\[C^{(\mathbf b)}(n,k)=p^{n-k}{n-1\choose k-1}.\]
\end{prp}
\begin{proof}
Equation (\ref{r1}) has the form:
\[C^{(\mathbf b)}(n,k)=\sum_{i=1}^{n-k+1}p^{i-1}C^{(\mathbf b)}(n-i,k-1).\;(k\leq n).\]
We prove the formula by induction on $k.$ It is obviously true for $k=1.$ Suppose it is also true for $k-1.$
Then the preceding equation takes the form:
\[C^{(\mathbf b)}(n,k)=p^{n-k}\sum_{i=1}^{n-k+1}{n-i-1\choose k-2}.\]
On the other hand, by a well-known horizontal recursion for the binomial coefficients we have
\[{n-1\choose k-1}=\sum_{i=1}^{n-k+1}{n-i-1\choose k-2},\]
and the formula is true.
\end{proof}
Using the binomial formula, for the number of all generalized compositions, we obtain
\[C^{(\mathbf b)}(n)=(1+p)^{n-1}.\]
This is the formula (i), Corollary  13, in \cite{mil}.

In the next two results we consider the case when $b_i$ is a
linear function of $i.$

\begin{prp}\label{prp2} Let $p,m,n$ be positive integers, and let $b_i=m(i-1),\;(i=1,2,\ldots).$ Then
 \[C^{(\mathbf b)}(n,k)=m^k\cdot{n-1\choose 2k-1}.\]
\end{prp}
\begin{proof} The proposition is obviously true for $k=1.$ Assume that it is true for $k-1.$
Equation (\ref{r1}) has the form:
\[C^{(\mathbf b)}(n,k)=m\sum_{i=1}^{n-k+1}(i-1)C^{(\mathbf b)}(n-i,k-1),\;(k\leq n).\]
Using the induction hypothesis yields
\[C^{(\mathbf b)}(n,k)=m^k\sum_{i=1}^{n-k+1}(i-1){n-i-1\choose 2k-3},\;(k\leq n).\]
It follows that
\[C^{(\mathbf b)}(n,k)=m^k\sum_{i=1}^{n-k}i{n-i-2\choose 2k-3},\;(k\leq n).\]

 Denote $S=\sum_{i=1}^{n-k}i{n-i-2\choose 2k-3}.$ Then,
 \[S=\sum_{i=0}^{n-k}{n-i-2\choose 2k-3}+\sum_{i=2}^{n-k}{n-i-2\choose 2k-3}+\cdots+\sum_{i=n-k}^{n-k}{n-i-2\choose 2k-3}.\]
Using the horizontal recursion for the binomial coefficients we obtain
\[S={n-2\choose 2k-2}+{n-3\choose 2k-2}+\cdots+{2k-2\choose
2k-2}.\] Using the same recursion once more we obtain
\[S={n-1\choose 2k-1}.\]
\end{proof}
For the number of all generalized composition we get
\begin{cor} In the conditions of Proposition \ref{prp2} we have
\[C^{(\mathbf b)}(n)=\sum_{k=1}^n{n-1\choose 2k-1}m^k.\]
\end{cor}

In a similar way we may prove the following:
\begin{prp} If $b_i=mi,\;(i=1,2,\ldots),$ then
 \[C^{(\mathbf b)}(n,k)=m^k\cdot{n+k-1\choose 2k-1}.\]
\end{prp}

Also,
\[C^{(\mathbf b)}(n)=\sum_{k=1}^n{n+k-1\choose 2k-1}m^k.\]
\begin{rem} Taking in particular $m=1$ in the preceding equation, we obtain Theorem 3.23, in \cite{hu}, about the so called $n$-colored compositions, defined in \cite{ag}.
\end{rem}
\section{Binomial coefficients}
In this section we investigate the generalized compositions, when the $b$'s are some binomial coefficients.
We first derive two closed formulas.
\begin{prp}\label{pb1} Let $k,p,n$ be positive integers, and let $b_i={p\choose i-1},\;(i=1,2,\ldots).$
Then,
 \[C^{(\mathbf b)}(n,k)={pk\choose n-k}.\]
 Also,
 \[C^{(\mathbf b)}(n)=\sum_{k=1}^n{pk\choose n-k}.\]
 \end{prp}
 \begin{proof} We go by induction on $k.$ For $k=1$ the proposition is obviously true. Using the induction hypothesis we see that the  first assertion is equivalent to the following identity:
\[ {pk\choose n-k}=\sum_{i=1}^{n-k+1}{p\choose i-1}{pk-p\choose n-i-k+1},\]
which is merely  the Vandermonde convolution.
 \end{proof}

The next result concerns the figured  numbers.
\begin{prp} Let $p,k,n$ be positive integers, and let \[b_i={p+i-1\choose p},\;(i=1,2,\ldots).\]
Then,
 \[C^{(\mathbf b)}(n,k)={n+pk-1\choose pk+k-1}.\]
 \end{prp}
 \begin{proof} We use induction on $k.$ For $k=1$ the proposition is obviously true. Using the induction hypothesis we see that the the assertion is equivalent to the following identity:
\[ {n+pk-1\choose pk+k-1}=\sum_{i=1}^{n-k+1}{p+i-1\choose p}{n-i+pk-p-1\choose pk-p+k-2}.\]
To prove this identity, we shall count $pk+k-1$-subsets of the set $X=\{1,2,\ldots,n+pk-1\}$ according to the place of its $(p+1)$the element in such a subset.
Suppose that this element is the $(p+i)$th element of $X.$ Such a subset may be chosen in ${p+i-1\choose p}\cdot{n-i+pk-p-1\choose pk-p+k-2}$ ways. We also conclude that $i$ ranges from $1$ to $n-k+1,$ which proves the proposition.

 \end{proof}

The following two results concern  the number of all generalized compositions.
We first prove that, in the case $b_{i}={i+p-1\choose
q},\;(i=1,2,\ldots),$ where $p,q$ are  positive integers, the  numbers $C^{(\mathbf b)}(n)$ satisfy a
homogenous linear recurrence  equation  of the $(q+1)$th order, with
constant coefficients.

\begin{prp}\label{kt} Let $p,q,n$ be  positive integers, and let $b_i={i+p-1\choose q},\;(i=1,2,\ldots).$
Then there exist integers $m_i(p,q),(i=0,1,\ldots,q),$ not depending on $n,$   such that
\begin{equation}\label{bi0}C^{(\mathbf b)}(n+q+1)=\sum_{i=0}^{q}m_i(p,q)C^{(\mathbf b)}(n+i),\;(n\geq 2).\end{equation}
\end{prp}
 \begin{proof}

We define the function $F(n,j)$ in the following way:
\begin{equation}\label{df}F(n,j)=\sum_{i=1}^{n-1}{n-i+p\choose q-j}C^{(\mathbf b)}(i-1),\end{equation}
where $0\leq j\leq q,\;2\leq n.$
We want to prove that the following equation holds
\begin{equation}\label{jf} F(n,j)=\sum_{i=0}^{j+1}c(i,j)C^{(\mathbf b)}(n+i-1),\end{equation}
where $c(i,j)$ are integers, depending only on $p$ and $q.$

The proof goes by induction on $j.$
Taking $n=1$ in (\ref{sve}) we get  $C^{(\mathbf b)}(1)={p\choose q}.$ For $n>1$
we get
\begin{equation}\label{bi1}C^{(\mathbf b)}(n)={p\choose q}C^{(\mathbf b)}(n-1)+\sum_{i=1}^{n-1}{n-i+p\choose q}C^{(\mathbf b)}(i-1).\end{equation}
It follows that
\begin{equation}\label{fn0}F(n,0)=C^{(\mathbf b)}(n)-{p\choose q}C^{(\mathbf b)}(n-1).\end{equation}
Hence, taking
\[c(0,0)=-{p\choose q},\;c(1,0)=1,\] we see that (\ref{jf}) holds for $j=0$ and $n\geq 2.$

Suppose that (\ref{jf}) holds for some $j\geq 0.$
Replacing $n$ by $n+1$ in (\ref{df}) yields
\[F(n+1,j)=\sum_{i=1}^{n}{n+1-i+p\choose q-j}C^{(\mathbf b)}(i-1).\]
Using the standard recursion for the binomial coefficients one obtains
\[F(n,j+1)=F(n+1,j)-F(n,j)-{p+1\choose q-j}C^{(\mathbf b)}(n-1).\]

Using the induction hypothesis yields
\[F(n,j+1)=\]\[=\sum_{i=0}^{j+1}c(i,j)C^{(\mathbf b)}(n+i)-\sum_{i=0}^{j+1}c(i,j)C^{(\mathbf b)}(n+i-1)-{p+1\choose q-j}C^{(\mathbf b)}(n-1).\]

Denoting
\[c(0,j+1)=-c(0,j)-{p+1\choose q-j},\;c(j+2,j+1)=c(j+1,j),\]
\[c(i,j+1)=c(i-1,j)-c(i,j),\;(1\leq i\leq j+1),\]
implies
\[F(n,j+1)=\sum_{i=0}^{j+2}c(i,j+1)C^{(\mathbf b)}(n+i-1),\;(n\geq 2),\]
and (\ref{jf}) is true.

Since $F(n,q)=\sum_{i=1}^{n-1}C^{(\mathbf b)}(i-1),$ we have
\begin{equation}\label{bi3}\sum_{i=0}^{q+1}c(i,q)C^{(\mathbf b)}(n+i-1)=\sum_{i=1}^{n-1}C^{(\mathbf b)}(i-1).\end{equation}
Replacing $n$ by $n+1$ in (\ref{bi3}) yields
\begin{equation}\label{bi4}\sum_{i=0}^{q+1}c(i,q)C^{(\mathbf b)}(n+i)=\sum_{i=1}^{n}C^{(\mathbf b)}(i-1).\end{equation}
Subtracting (\ref{bi4}) from (\ref{bi3}) we obtain
\begin{equation}\label{bi5}\sum_{i=0}^{q+1}c(i,q)\left[C^{(\mathbf b)}(n+i-1)-C^{(\mathbf b)}(n+i)\right]+C^{(\mathbf b)}(n-1)=0.\end{equation}

Further, we obviously  have $c(q+1,q)=1.$ Also, we may easily obtain the
values for $c(0,q+1).$ First, we have
\[c(0,1)=-c(0,0)-{p+1\choose q}={p\choose q}-{p+1\choose q}=-{p\choose q-1}.\]
Using  induction easily implies that
\begin{equation}\label{0j}c(0,j)=-{p\choose q-j},\;(j=0,1,\ldots,q).\end{equation}
In particular, $c(0,q)=-1,$ which means that $C^{(\mathbf b)}(n-1)$ vanishes in
equation (\ref{bi5}). Hence, equation (\ref{bi5}) becomes
(\ref{bi0}), if we take
\[m_i(p,q)=-c(i+1,q+1),\;(i=0,1,\ldots,q).\]
\end{proof}

\begin{rem} We have seen, in Proposition \ref{pxx}, that in the case $p-1=q,$ the number $C^{(\mathbf b)}(n)$ is the number of $q$-matrix compositions, as they are defined in \cite{mu}. Thus the numbers of $q$-matrix compositions satisfy a $(q+1)$th order homogenous linear recurrence equation with constant coefficients.
\end{rem}

\begin{rem} The coefficients $c(i,j),\;(j=0,1,\ldots;\;i=0,1,\ldots,j+1)$ form a kind of a Pascal-like triangle.
\end{rem}

We shall now consider the particular case $p=1,\;q>1,$ and  show  that then  the coefficients $m_i(1,q)$ can be obtained  explicitly.

\begin{prp}\label{kt1} Let $q$ be a positive integer, and let $b_i={i\choose q},\;(i=1,2,\ldots).$
Then,
\[C^{(\mathbf b)}(n+q+1)=\sum_{i=0}^{q}(-1)^{i+q}{q+1\choose i}C^{(\mathbf b)}(n+i)+
C^{(\mathbf b)}(n+1),\;(n\geq 2).\]
\end{prp}
\begin{proof}

Firstly, we have
\[c(0,0)=0,\;c(1,0)=1.\]
For $j\geq 1$, by (\ref{0j}), we have
\[c(0,j)=-{1\choose q-j}.\]
 It follows that
\[c(0,q-1)=c(0,q)=-1,\mbox{ and }c(0,j)=0\mbox{ otherwise }.\]
Furthermore, for $j<q$ we have
\[c(1,j)=c(0,j-1)-c(1,j-1)=-c(1,j-1)=c(1,j-2)=\ldots=(-1)^j,\]
and
\[c(1,q)=c(0,q-1)-c(1,q-1)=-1-c(1,q-1)=\ldots=-1+(-1)^q.\]
Also,
\[c(2,j)=(-1)^{j-1}j,\;(j\leq q).\]

We next prove that for $j,$ satisfying the condition $2\leq j\leq q,$ we have
\[c(i,j)=(-1)^{j-i+1}{j\choose i-1},\;(i=2,\ldots,j).\]
The equation is true for $i=2,$ by the preceding equation.
Suppose that it is true for some $i-1\geq 2.$ From the equation
\[c(i,j)=c(i-1,j-1)-c(i,j-1),\]
using the induction hypothesis we obtain
\[c(i,j)=(-1)^{j-i+1}{j-1\choose i-2}-c(i,j-1).\]

 From this we easily conclude that
 \[c(i,j)=(-1)^{j-i+1}\left[{j-1\choose i-2}+{j-2\choose i-2}+\cdots+{i-2\choose i-2}\right].\]
 The assertion is true, by the horizontal recursion for the binomial coefficients.
In particular, we have
\[(-1)^{i+q}[c(i+1,q)-c(i,q)]={q\choose i}+{q\choose i-1}={q+1\choose i}.\]
\end{proof}
Now, we shall derive the closed formula for the recursion from the preceding proposition.

\begin{prp} Let $q$ be a positive integer, and  let $b_i={i\choose q},\;(i=1,2,\ldots).$ Then,
 \[C^{(\mathbf b)}(n,k)={n+k-1\choose qk+k-1}.\]
\end{prp}
\begin{proof} We first conclude that each term of any generalized composition is $\geq q.$ It follows that  $C^{(\mathbf b)}(n,k)=0,$ if $n<qk.$ This means that the assertion holds for $n<qk.$  Assume that $n\geq qk.$

Using induction we easily conclude that the  assertion is equivalent to the following binomial identity:
\[{n+k-1\choose qk+k-1}=\sum_{i=1}^{n-k+1}{i\choose q}{n+k-2-i\choose qk-q+k-2},\;(qk\leq n).\]
Adjusting the lower and the upper bounds in the sum on the right-hand side,  we obtain the following identity:
\[{n+k-1\choose qk+k-1}=\sum_{i=q}^{n-qk+q}{i\choose q}{n+k-2-i\choose qk-q+k-2},\;(qk\leq n).\]

To prove this identity we shall count $(qk+k-1)$- subsets of the set $X=\{1,2,\ldots,n+k-1\}$
in the following way:
Suppose that  $x$ is the $(q+1)$th element of a $(qk+k-1)$-subset of $X,$ and suppose that  we have $i$ elements of $X$ in the subset,   which are less than $x.$ It follows that  there are \[{i\choose q}{n+k-2-i\choose qk-q+k-2}\] subsets with this property.
The assertion is true, since $i$ ranges from $q$ to $n-qk+q.$
\end{proof}
As an immediate consequence we have
\begin{cor} If $b_i={i\choose q}, (i=1,2,\ldots),$ then
\[C^{(\mathbf b)}(n)=\sum_{k=1}^n{n+k-1\choose qk+k-1}.\]
\end{cor}
\begin{rem} The preceding equation is the closed formula for the recurrence equation from Proposition \ref{kt1}.
\end{rem}

\section{Catalan numbers}
In this section we consider the case when the $b$'s are Catalan numbers.
In the first result we shall prove that the numbers of generalized compositions with a fixed number of parts, may be expressed
in terms of the numbers of the so called Catalan triangle,
introduced by Chapiro, \cite{ch}. We let $\mathbf c_i$  denote the $i$th Catalan number. Also, $B(n,k)$ denotes a number of Catalan triangle.
Thus,
\[B(n,k)=\frac{k}{n}{2n\choose n+k},\;(k\leq n).\]

\begin{prp}\label{k1} Let $n,k$ be positive integers,  and let $b_i=\mathbf c_{i},(i=1,2,\ldots).$ Then,
\[C^{(\mathbf b)}(n,k)=B(n,k).\]
Further,
\[C^{(\mathbf b)}(n)={2n-1\choose n}.\]
\end{prp}
\begin{proof}
Equation (\ref{r1}), in this case,  has the form:
\[C^{(\mathbf b)}(n,k)=\sum_{i=1}^{n-k+1}\mathbf c_{i}C^{(\mathbf b)}(n-i,k-1),\;(k\leq n).\]
The assertion follows by induction, using Theorem 14.3, \cite{ks}.
The second assertion follows from Theorem 14.2, \cite{ks}.
\begin{rem}
Note that, in the preceding proposition, we have an example when the number of all generalized compositions is a binomial coefficient.
\end{rem}
\end{proof}
We now slightly change the conditions of the preceding corollary to obtain a relationship among Catalan numbers, binomial coefficients, and the numbers of Catalan triangle.

 \begin{prp} Let $n,k$ be positive integers, and let $b_i=\mathbf c_{i},\;(i=0,1,\ldots).$ Then, for $n\geq k,$ we have
\begin{equation}\label{kb}C^{(\mathbf b)}(n,k)=\sum_{i=0}^{k-1}{k\choose i}B(n-k,k-i).\end{equation}
\end{prp}
\begin{proof}
We shall first prove that, for $1\leq k\leq n,$ the following
equation
\begin{equation}\label{xx}C^{(\mathbf
b)}(n,k)=\sum_{i_1+i_2+\cdots+i_k=n-k}\mathbf c_{i_1}\cdot \mathbf
c_{i_2}\cdots \mathbf c_{i_k},\end{equation}  holds. The sum is
taken over $i_1\geq 0,i_2\geq 0,\ldots,i_k\geq 0.$
 We use induction on $k.$
For $k=1,$ by (\ref{r1}), we have $C^{(\mathbf b)}(n,1)=\mathbf
c_{n-1}.$ On the other hand, (\ref{xx}) has the form: \[
C^{(\mathbf b)}(n,1)=\sum_{i_1=n-1}c_{i_1}=\mathbf c_{n-1},\] and
the proposition is true. Suppose that the proposition is true for
$k\geq 1.$ Then,
 \[C^{(\mathbf b)}(n,k+1)=\sum_{i=1}^{n-k}\mathbf c_{i-1}C^{(\mathbf b)}(n-i,k).\]

 Using the induction hypothesis yields
\[C^{(\mathbf b)}(n,k+1)=\sum_{i=1}^{n-k}\mathbf c_{i-1}\sum_{i_1+i_2+\cdots+i_k=n-i-k}\mathbf c_{i_1}\cdot \mathbf c_{i_2}\cdots \mathbf c_{i_k}.\] Denote $i-1=i_{k+1}$ to obtain
\[C^{(\mathbf b)}(n,k+1)=\sum_{i_1+i_2+\cdots+i_k+i_{k+1}=n-k-1}\mathbf c_{i_1}\cdot \mathbf c_{i_2}\cdots \mathbf c_{i_k}\cdot \mathbf c_{i_{k+1}},\]
and (\ref{xx}) is true.

 Collecting terms with a fixed number of zeroes in (\ref{xx}) we obtain
 \[C^{(\mathbf b)}(n,k)=\sum_{j=0}^{k-1}{k\choose j}\sum_{i_1+i_2+\cdots+i_{k-j}=n-k}\mathbf c_{i_1}\cdot \mathbf
c_{i_2}\cdots \mathbf c_{i_{k-j}},\] where all sums on the
right-hand side are taken over $i_t\geq 1.$ According to Theorem 14.
4, \cite{ks},  we have

\[B(n,k)=\sum_{i_1+i_2+\cdots+i_{k}=n}\mathbf c_{i_1}\cdot \mathbf
c_{i_2}\cdots \mathbf c_{i_{k}},\]
 where $i_1\geq 1,\ldots,i_k\geq 1,$ and the proposition is true.
\end{proof}
In \cite{mil} it is proved that the sum on the right-hand side of
equation (\ref{xx})  equals the number of the weak compositions of
$n-k$ in which exactly $k$ parts equal $0.$  We thus have
\begin{cor} Let $n,k$ be positive integers, and let $b_i=\mathbf c_{i},\;(i=0,1,\ldots).$ Then
 $C^{(\mathbf b)}(n,k)$ is the number of the weak generalized
compositions of $n-k$ in which there are exactly $k$ zeroes.
\end{cor}

It is proved in Proposition 3, \cite{mil}, that in this case $\mathbf c_n$ is the number of all generalized compositions. We thus obtain a formula
which shows that Catalan numbers are some kind of convolution of the numbers
of Pascal and Catalan triangles.

\begin{cor} Let $n$ be a positive integer. Then
\[\mathbf c_{n}=1+\sum_{k=1}^{n-1}\sum_{i=1}^{k-1}{k\choose i}B(n-k,k-i).\]
\end{cor}

\end{document}